

\baselineskip=14pt
\parskip=10pt
\def\halmos{\hbox{\vrule height0.15cm width0.01cm\vbox{\hrule height
  0.01cm width0.2cm \vskip0.15cm \hrule height 0.01cm width0.2cm}\vrule
  height0.15cm width 0.01cm}}
\font\eightrm=cmr8 
\font\eighttt=cmtt8
\magnification=\magstephalf

\def\1{{\overline{1}}}
\def\2{{\overline{2}}}
\parindent=0pt
\overfullrule=0in
\def\frac#1#2{{#1 \over #2}}
\bf
\centerline
{
A Proof of the Loehr-Warrington
Amazing TEN to the Power n Conjecture
}
\rm
\bigskip
\centerline{ {\it
Shalosh B. EKHAD${}^1$, Vince VATTER${}^1$, and Doron 
ZEILBERGER}\footnote{$^1$}
{\eightrm  \raggedright
Department of Mathematics, Rutgers University (New Brunswick),
Hill Center-Busch Campus, 110 Frelinghuysen Rd., Piscataway,
NJ 08854-8019, USA.
{\eighttt [ekhad, vatter, zeilberg] at math dot rutgers dot edu} ,
\hfill \break
{\eighttt http://www.math.rutgers.edu/\~{}[zeilberg/ekhad.html, 
vatter, zeilberg]} .
First version: Sept. 14, 2005.
Accompanied by Maple package {\eighttt TEN}
downloadable from Zeilberger's website.
Supported in part by the NSF.
}
}

{\bf Theorem:}
There are $10^n$ words in the alphabet $\{3,-2\}$ of length
$5n$, sum $0$, and such that
every factor that sums to $0$ and that starts with a $3$ may not
be immediately followed by a $-2$.

{\bf Proof:}
{\bf 1.} Download Maple package {\tt TEN} from
{\eighttt http://www.math.rutgers.edu/\~{}zeilberg/tokhniot/TEN}.
Save it as {\tt TEN}.
{\bf 2.} Go into Maple, and type: {\tt read TEN:}.
{\bf 3.} Let your computer rediscover a (linear) grammar
for these words, by typing
{\tt G:=DiscoverGrammar(3,2,3):}
(if you actually want to see it replace the : by ; ).
{\bf 4.} To prove (rigorously!) that the conjectured grammar $G$ indeed 
describes
(unambiguously!) our language, type:
{\tt ProveGrammar32(G,3); } and get the output {\tt true}.
{\bf 5.} To get the weight-enumerator for the language (where the
weight of a word $w$ is $x^{length(w)}$), type
{\tt GFgrammar(G,x); }
and get the output $1/(1-10x^5)$. \halmos\ ({\bf endproof})
\halmos\ ({\bf endpaper})

{\bf Appendix: TEN Remarks for Human Readers}

{\bf 1.} The Theorem was conjectured by Nick Loehr and Greg Warrington, 
in a more general setting,
where $3$ and $-2$ are replaced by general relatively-prime positive
and negative integers, $a$ and $-b$, $5n$ is replaced by $(a+b)n$,
and $10^n$ is replaced by ${{a+b} \choose {a}}^n$. The general
case is still open. Greatly inspired by our proof,
Loehr and Warrington, together with Bruce Sagan,
found a computer-free proof of our theorem. They then extended
their approach to prove it for $b=2$ and all odd $a$.
This will appear in their forthcoming paper
{\it ``A Human Proof for a Generalization of Shalosh B. Ekhad's $10^n$
Lattice Paths Theorem.''}
We wish to acknowledge stimulating discussions with these three humans,
and of course we thank Nick and Greg for conjecturing such a beautiful
result in the first place, and Bruce for telling
Vince, who told Doron, who told Shalosh.

{\bf 2.} The Maple package {\tt TEN} that
(automatically!!)~discovered the grammar, and then
(automatically!!!)~proved its correctness, and then
(automatically!)~computed the weight-enumerator, was
written by VV and DZ. It was executed by SBE.

{\bf 3.} For a blow-by-blow description of how the
grammar was {\it discovered}, type

{\tt G:=DiscoverGrammarVerbose(3,2,3):} .

(The output file may be viewed/downloaded at
{\eighttt http://www.math.rutgers.edu/\~{}zeilberg/tokhniot/oTENdgv}.)

{\bf 4.} People with no access to Maple may view the output
at  \hfill\break
{\eighttt http://www.math.rutgers.edu/\~{}zeilberg/tokhniot/oTEN32t},  \hfill\break
{\eighttt http://www.math.rutgers.edu/\~{}zeilberg/tokhniot/oTEN32v},  \hfill\break
and {\eighttt http://www.math.rutgers.edu/\~{}zeilberg/tokhniot/oTEN32vv}
for the terse, verbose, and very verbose versions.

{\bf 5.} The best way to understand the {\it heuristics} behind the 
{\bf act of discovery}
(of the linear grammar for the language of the theorem)
and the {\it logic} behind the {\bf act of verification} (proving 
rigorously
that the empirically-conjectured grammar is indeed correct) is to 
carefully read  the
{\bf Maple source-code}, generously made available by the authors, free 
of charge,
even though the same {\bf methodology} is very likely to solve other 
problems.
There are ample comments, and readers are encouraged to experiment with 
{\tt TEN}
themselves. A nice exercise would be to type {\tt 
DiscoverGrammar(1,1,3);}
for the grammar of the language of zero-sum words in $\{-1,1\}$ that 
avoid
factors of the form $1[-1](-1)$, where $[-1]$ denotes a word that sums 
up to $-1$.
Then prove its validity by human means.

{\bf 6.} What if you don't know Maple well enough to follow someone 
else's code?
Then go and learn Maple! It would be time much better spent than reading
esoteric papers like this one.

{\bf 7.} Having said that, as a concession to Maple-illiterate people,
let us briefly describe the methodology of {\bf discovery}, to be
followed by a description of the  methodology of {\bf proof}.

Procedure {\tt Corpus} first  generates all the words of our language 
up to
a specified length. For any language $L$ in the alphabet $\{a,-b\}$, 
and any
pair of words in that alphabet (not necessarily in $L$),
$[w_1,w_2]$, let $L(w_1,w_2)$ be the subset of $L$ consisting
of those words of the form $w_1uw_2$ for some word $u$.
We want to construct a binary {\it family tree}, rooted at 
$L=L(\phi,\phi)$,
whose vertices are
pairs $[w_1,w_2]$ that stand for $L(w_1,w_2)$.
It may happen that $L(w_1,w_2)$ is {\bf empty}, i.e.~$L$ has
no words of the form $w_1uw_2$.
It may also happen that such an $L(w_1,w_2)$ is a {\bf clone} of
another $L(w'_1,w'_2)$ , i.e.
$$
\{ v \, \vert \, w_1 v w_2 \in L \}=\{ v \, \vert \,  w'_1 v w'_2 \in L 
\} \quad ,
$$
where $length(w'_1)+length(w'_2)<length(w_1)+length(w_2)$.

Each $L(w_1,w_2)$ may be naturally partitioned in two different ways.
$$
L(w_1,w_2)=L(w_1a,w_2) \cup L(w_1(-b),w_2) \quad (\cup \{w_1w_2\}, \,\, 
if \,\,  w_1w_2 \in L ) \quad ,
\eqno(HeadWay)
$$
since the letter immediately following $w_1$ is either $a$ or $-b$
(assuming that the middle part is not null, in which case we have to add
the additional singleton)
or
$$
L(w_1,w_2)=L(w_1,aw_2) \cup L(w_1,(-b)w_2) \quad (\cup \{w_1w_2\},\,\, 
if \,\,w_1w_2 \in L )\quad,
\eqno(TailWay)
$$
since the letter immediately preceding the  $w_2$ is either $a$ or $-b$
(assuming that the middle part is not null, in which case we have to add
the additional singleton).

Anthropomorphizing a bit, we can think of every vertex as a man who has
two wives, let's call them Rachel and Leah. Unlike Jacob,
this man doesn't decide whom he likes better until he sees the sons 
that they
give him. Rachel gives birth to Joseph and Benjamin, and Leah gives 
birth
to Reuben and Simon (let's pretend that Levi et al. never got born).
Now, deciding by the children's merit, he picks one of the wives as
the main one and the other one becomes his concubine, and her sons
get disowned. He may not mix sons!  He either picks
Joseph and Benjamin as heirs, or Reuben and Simon, but he is not
allowed to pick, say, Joseph and Simon.

Having picked Joseph and Benjamin as Jacob's legal heirs, they
are now leaves in our expanding tree.
Consider such a  new vertex (that starts out as a bachelor).
It may be empty,
or it may be a {\bf clone} of some older (legitimate!) relative (not 
necessarily
a direct ancestor, e.g. Abraham is definitely okay, and Isaac, of 
course,
but Essau would also do).
In that case it becomes a {\bf permanent leaf}
(confirmed old bachelor),  and is {\bf forbidden} to have sons.
But if it is neither empty nor a clone, then he, in his turn,
has two son-pairs and he
must decide which pair to pick as legal heirs,
in other words, he has to make up his mind whether
to split according to $(HeadWay)$ or according to $(TailWay)$.

The way procedure {\tt DiscoverGrammar} in our Maple package {\tt TEN} 
decides
this issue is by preferring those sons
with congenial cardinalities, i.e.~whose
greatest-common-divisor with powers of ${{a+b} \choose {b}}$
is as large as possible. We admit that it is only
one possible heuristics for picking heirs, and its only merit is that 
it worked.
The process terminates when all the leaves are either empty or clones.

Having discovered the grammar, to
find the {\bf weight-enumerator} of its language,
(this is implemented in procedure {\tt GFgrammar} of {\tt TEN}),
we set up a system of linear equations, whose unknowns are
the weight-enumerators of the $L(w_1,w_2)$'s, let's call them
$z[w_1,w_2]$, for all vertices $[w_1,w_2]$ (both internal vertices and 
leaves).
For each internal vertex $[w_1,w_2]$, we have the equation
$$
z[w_1,w_2]=z[w_1',w_2']+z[w_1'',w_2''] \quad ,
\eqno(Internal)
$$
where $[w_1',w_2']$, $[w_1'',w_2'']$ are the two (legitimate!) children 
of $[w_1,w_2]$.
If $w_1w_2 \in L$, we have to use, instead
$$
z[w_1,w_2]=z[w_1',w_2']+z[w_1'',w_2'']+x^{len(w_1)+len(w_2)} \quad .
\eqno(Internal')
$$
For each leaf, $[w_1,w_2]$, if it is empty, we have the obvious equation
$$
z[w_1,w_2]=0 \quad ,
\eqno(LeafEmpty)
$$
while if it is a clone of $[w_1',w_2']$, say, then we have
$$
z[w_1,w_2]=x^{len(w_1)+len(w_2)-len(w_1')-len(w_2')}z[w_1',w_2'] \quad.
\eqno(LeafClone)
$$
Now we (or rather Maple) solve(s) this huge system of equations and
get(s) the weight-enumerators of all vertices, in particular,
$z[\phi,\phi]$, our object of desire.

{\bf 8.} We still need to {\bf rigorously} prove that the language 
generated by our
putative grammar is indeed the language of interest. Since  the
paternity part is obviously true, we have to prove leafness,
both of the empty and the clone kinds.
To prove that a supposedly empty leaf is indeed so, we must
demonstrate that for each such (supposedly)
empty leaf $[w_1,w_2]$, there can never be a {\it good} word (i.e.~a 
word in $L$, of {\it whatever} length),
of the form $w_1vw_2$.  In other words we have to
{\bf logically} prove the implication:

Every word of the form $w_1 v w_2$ must be bad, i.e.~contain a mishap, 
that is  a factor of the form
$a[-a](-b)$.  (Here, for any integer $A$, $[A]$ denotes any word (in 
the alphabet)
that adds up to $A$.)

Note that if $w_1 v w_2$ belongs to $L$ then the sum of $v$ is $-A$,
where $A:=sum(w_1)+sum(w_2)$, so we have to prove that
every word of the form $w_1[-A]w_2$ {\it must contain a mishap},
i.e.~a factor of the form $a[-a](-b)$.

The claim that vertex $[w_1,w_2]$ is a clone of vertex $[w'_1,w'_2]$ is:

``$w_1vw_2$ is good iff $w'_1vw'_2$ is good'',
or equivalently,

``$w_1vw_2$ is bad iff $w'_1vw'_2$ is bad'',

which really contains two statements

``if $w_1vw_2$ is bad  then $w'_1vw'_2$ is bad'', and

``if $w'_1vw'_2$ is bad  then $w_1vw_2$ is bad''.

To prove such inclusions, the computer looks at all {\it potential 
mishaps},
which consist of an actual `$a$' in the $w_1$ part and a potential 
`$-b$' in the $v$ part,
or a potential `$a$' in the $v$ part and an actual `$-b$' in the $w_2$ 
part.
Also, conceivably (but rarely) we should consider an actual $a$ in the 
$w_1$
part and an actual $-b$ in the $w_2$ part.
Each such potential mishap entails some factorization of $v$ of the form
$[A]a[B]$ or $[A](-b)[B]$, or in the last case, just plain $[A]$.
This is implemented (for the $\{3,-2\}$ case) in procedure {\tt 
PotentialMishaps32}
of {\tt TEN}.

At the bottom line, proving that $L(w_1,w_2)$ is {\it empty} boils down
to proving that every word of the form $w_1[A]w_2$, where
$A:=-(sum(w_1)+sum(w_2))$, {\bf must} contain some mishap.
Also proving cloneship reduces to proving several statements
of the form: ``$u_1 [A] u_2 [B] u_3$ {\bf must} always contain some 
mishap'', where $u_1,u_2,u_3$ are
{\it specific words}, and $[A],[B]$ denote arbitrary words that sum to 
$A$ and $B$ respectively,
for some specific integers $A$ and $B$.
Of course, w.l.o.g.~the $[A]$ and $[B]$ words are mishap-free.

A crucial tool in the automated proof of such  assertions
is a recurrence that is a simple consequence of the
{\bf Discrete Rolle Theorem}:
$$
[A]=
\cases{
\displaystyle\bigcup_{i=1}^{a} \,\, [i-a]a[A-i] ,& if $A > 0$;\cr
\quad   &  \quad \cr
\displaystyle\bigcup_{i=1}^{b} \,\, [-i+b](-b)[A+i],
& if $A<0$.\cr}
\eqno(WordRecurrence)
$$
Indeed any word that sums to $A>0$ must have its {\bf shortest} prefix 
whose
sum is positive, the last letter of that prefix being, necessarily, $a$
(or else it wouldn't be the shortest), and the sum of that prefix
(by `continuity'!) must be between $1$ and $a$
(inclusive). Analogously, for $A<0$.

By iterating these recurrences it is easily seen
that each word $[A]$ can be written as one of several {\it 
factorizations}
featuring the letters $a$ and $-b$ and the {\bf fundamental factors} 
$[0],[1], \dots, [b-1]$.
We can even apply $(WordRecurrence)$ to these, but then we'll
get `self-referential' expressions, i.e.~ways of expressing them in 
terms of
themselves (and the actual letters $a$ and $-b$) as a kind of
`Chomskian' grammar.

Whenever we have to prove emptyness or cloneship, we have to prove
that a certain set of words doesn't exist. Like in 4CT or any
{\it reducto} proof (all the way back to $\sqrt{2} \not \in Q$), if the 
set
was not empty, there would be a {\bf minimal counterexample}.
So it suffices to prove that minimal counterexamples do not exist.
Note that a  minimal counterexample (in our business) can never have
a {\it proper} zero-factor (removing a zero-factor from a good word
obviously results in another good word).
So if we are in the lookout for minimal
counterexamples, we can considerably reduce the many  options that  an 
$[A]$ word can
have, by discarding those factorizations in the
enlarged (`meta-') alphabet $\{a,-b,[1],[2], \dots, [b-1]\}$ that 
contain a zero-sum proper factor.
Of course we also kick out any such meta-word that contains a mishap,
i.e.~a factor of the form $a[-a](-b)$.
If there are any survivors left after the first purge, we can replace
$[1], \dots, [b-1]$ by their self-referential expressions
getting longer `meta-words', and once again
discard all those that contain mishaps and all those that
contain a zero-sum proper factor. If all goes well, nothing
will be left after finitely many purges.
For the $\{3,-2\}$ case three purges sufficed to prove every
instance of emptyness or cloneship.

Let's explicate the above for our case of interest, $a=3$, $b=2$.
There is only one `meta-letter', [1], and we have, according to {\tt 
TEN}
(but this is so simple that even you can do it!)
$$
[1]=\{[3, -2], [-2, 3], [-2, [1], 3, [1], -2]\} \quad,
$$
which means that every word in the alphabet $\{3, -2\}$
that contains no zero-sum-factors and no mishaps
is either the two-letter word $(3)(-2)$ or the two-letter
word $(-2)3$ or else can be written as
$(-2)w_13w_2(-2)$, where $w_1,w_2$ are zero-sum-free and mishap-free 
words
whose sum is $1$ on their own right.

We also have:
$$
[-1]=\{[-2, [1]], [[1], -2]\} \quad,
$$
stating that any mishap-free and zero-factor-free
word in $\{3,-2\}$ that sums up  to $-1$ may be written
as $(-2)w_1$ or $w_1(-2)$ where $w_1$ is a such a good word
that adds up to $1$.
You are welcome to decipher the following lemmas that
Shalosh discovered and proved.
$$
[-2]=\{[-2], [[1], -2, -2, [1]]\} \quad ,
$$
$$
[2]=\{[3, -2, [1]], [3, [1], -2], [-2, 3, [1]], [-2, [1], 3], [[1], -2, 
3]\} \quad ,
$$
$$
[3]=\{[3], [-2, [1], 3, [1]], [[1], -2, 3, [1]], [-2, 3, 3, -2, [1]],
     [-2, 3, 3, [1], -2], [-2, 3, -2, 3, [1]]\} \quad .
$$

{\bf 9.} The grammar is especially simple for the case $b=1$.
If $A_a(n)$ is the set of good words in $\{a,-1\}$ we have
(here $\1$ stands for $-1$)
$$
A_a(n)=A_a(n-1) \, a\, \1^a \,\, \cup \,\,
\1\,A_a(n-1)\, a\, \1^{a-1} \, \cup \,
\1\,\1A_a(n-1)a\,\1^{a-2} \, \cup \, \dots \, \cup \,
\1^a\,A_a(n-1)\,a \quad ,
$$
immediately implying that $\vert A_a(n) \vert =(a+1)^n$.

{\bf 10.} ``{\bf Mr. Watson, come here. I want you!}''

While grammatically correct, this sentence is not quite Shakespeare.
Its great significance is the {\bf way} it was uttered, via
the telephone, by Alexander Graham Bell. Analogously, while our proof
is definitely closer to Gauss than Bell's sentence is to Shakespeare,
its main significance, if we do say so ourselves, is in the
{\bf way} that it was discovered, and especially the {\bf way} that
it was {\bf proved}, all by computer! The role of the humans (VV and DZ)
is no longer that of {\it athlete}, but that of {\it coach}, no longer
{\it prover} but {\it programmer}, or if you wish, {\it meta-prover}.
We believe that this will be {\bf the  way to go} sooner than you 
think! Amen.

\end